\documentclass[12pt,a4paper]{article}
\usepackage{amssymb,amsmath}
\usepackage{emlines,bezier}
\newtheorem{theorem}{\sl THEOREM }[section]

\newtheorem{remark}{\sl REMARK }[section]
\newtheorem{corollary}{\sl COROLLARY }[section]
\newtheorem{defini}{\sl DEFINITION}

\renewcommand{\S}{\mathhexbox{278}\ }
\newcommand{\epr}{\hfill $\square$}
\newcommand{\R}{\mathbb R}
\newcommand{\G}{\mathcal G}
\newcommand{\Sc}{\mathcal S}
\newcommand{\Lc}{\mathcal L}
\newcommand{\Fc}{\mathcal F}
\begin{document}
\begin{center}
\bf ON THE NONHAMILTONIAN CHARACTER OF SHOCKS
IN 2-D PRESSURELESS GAS
\end{center}
\vskip 6 pt
\centerline{Yu.G.Rykov}
\vskip 12pt

%
%  ITALIAN ABSTRACT
%
{\scriptsize {\sl PRESENTAZIONE:}\it\ Si considera un sistema
bidimensionale della dinamica dei gas introdotto nel 1970 da Ya.~Zeldovich
per descrivere la formazione della struttura di grande scala dell'
universo.  Il sistema si rivela come qualcosa di intermedio tra un sistema
di equazioni differenziali ordinarie e un sistema iperbolico di equazioni
alle derivate parziali. La caratteristica principale e' la nascita di
singolarita': discontinuita' della velocita' e funzioni $delta$ di vario
tipo per la densita'. Si da' una descrizione rigorosa delle soluzioni
generalzzate in termini di misure di Radon e si ottiene una
generalizzazione delle condizione di Rankine-Hugoniot.
Sulla base di tali condizioni si mostra che la rappresentazione
variazionale delle soluzioni generalizzate, valida nel caso
unidimensionale, non vale in generale nel caso bidimensionale. Si ottiene
anche un sistema unidimensionale non banale non strettamente iperbolico
per la descrizione dell' evoluzione delle {\rm all' interno} dell' urto.}

%
%  ABSTRACT
%
{\scriptsize {\sl ABSTRACT:}\it\ The paper deals with the 2-D system of gas
dynamics without pressure which was introduced in 1970 by Ya.~Zeldovich to
describe the formation of large-scale structure of the Universe. Such
system occurs to be an intermediate object between the systems of
ordinary differential equations and hyperbolic systems of PDE. The main
its feature is the arising of singularities: discontinuities for velocity
and $\delta$-functions of various types for density. The rigorous notion
of generalized solutions in terms of Radon measures is introduced and the
generalization of Rankine-Hugoniot conditions is obtained. On the basis of
such conditions it is shown that the variational representation for the
generalized solutions, which is valid for 1-D case, in 2-D case generally
speaking does not take place. A nontrivial 1-D system of nonstrictly
hyperbolic type is also obtained to describe the evolution {\rm inside}
the shock.}

\tableofcontents
\allowdisplaybreaks
%
%  SECTION 1
%
\section{Introduction. Basic definitions.}

This paper studies the shock waves for the system of 2-D pressureless
gas dynamics
\begin{equation} \left\{ \begin{aligned}
\ &\frac{\partial\varrho}{\partial t}+\frac{\partial (\varrho u)}
{\partial x}+\frac{\partial (\varrho v)}{\partial y}=0 \\
\ &\frac{\partial (\varrho u)}{\partial t}+\frac{\partial (\varrho u^2)}
{\partial x}+\frac{\partial (\varrho u v)}{\partial y}=0 \\
\ &\frac{\partial (\varrho v)}{\partial t}+\frac{\partial (\varrho u v)}
{\partial x}+\frac{\partial (\varrho v^2)}{\partial y}=0,
\qquad (t,x,y) \in \R_+ \times \R^2\ ,
\end{aligned} \right. \label{1} \end{equation}
where $(u,v)$ has the physical meaning of velocity vector, $\varrho$ has
the physical meaning of density. The system (\ref{1}) is nonstrictly
hyperbolic system of conservation laws which has three coinciding
characteristic fields and incomplete set of eigenvectors. Due to these
properties the shock waves develop strong singularities in the density
which are of type of $\delta$-functions on the surface. So (\ref{1}) is
relevant to describe some process of concentration of matter. We will
consider the Cauchy problem for (\ref{1})
\begin{align}
\varrho(0,x,y) &= \varrho_0(x,y)>0 \notag \\
u(0,x,y) &= u_0(x,y) \label{2} \\
v(0,x,y) &= v_0(x,y)\ , \notag
\end{align}
where $\varrho_0,u_0,v_0$ are piecewise $C^1(\R^2)$ functions and will be
taken in the special form (see below) to produce locally the single shock
front.

For the smooth functions the system (\ref{1}) is equivalent to the following
system of equations
\begin{equation} \left\{ \begin{aligned}
\ &\frac{\partial\varrho}{\partial t}+\frac{\partial (\varrho u)}
{\partial x}+\frac{\partial (\varrho v)}{\partial y}=0 \\
\ &\frac{\partial u}{\partial t}+u \frac{\partial u}{\partial x}
+v \frac{\partial u}{\partial y}=0 \\
\ &\frac{\partial v}{\partial t}+u \frac{\partial v}{\partial x}
+v \frac{\partial v}{\partial y}=0, \qquad (t,x,y)
\in \R_+ \times \R^2\ .
\end{aligned} \right. \label{3} \end{equation}

The last two equations in the system (\ref{3}) constitute the so-called
inviscid 2-D Burgers equation, which was proposed (but in 3-D case) by
Ya.~Zeldovich \cite{ZeldAsAs70} to describe the formation of the large
scale structure of the Universe. Further this approach, in particular the
consideration of the whole system (\ref{3}), was developed in the
consequent papers (see, for example, \cite{ShZeRMPh89}, \cite{GMSaBNon91},
\cite{SAFrCMPh92}, \cite{VDFNAsAs94} and the references therein) from the
physical point of view. But if we are dealing with the laws of
conservation of mass and momentum it would be more convenient from the
mathematical and physical also points of view to investigate the system in
divergent form (\ref{1}) especially in the case when one has the
developing of such singularities as shock waves. For 1-D variant of the
system (\ref{1}) the particular generalized solutions were constructed in
\cite{BoucAKiC94} (there also were proposed some numerical schemes), and
the existence theorem in the sense of {\it Radon measures} for wide class
of initial data was proved in \cite{ERSiUMNa95}, \cite{ERSiCMPh96},
\cite{GrenCRAS95}, see also \cite{BouJCRAS98}, \cite{BouJCPDE99} for the
uniqueness results in the framework of {\it duality} solutions. In
\cite{ERSiUMNa95}, \cite{ERSiCMPh96} the generalized solutions were
constructed by the aid of variational principle for appropriate
Hamilton-Jacobi equation (see \cite{HopfCPAM50}, \cite{Lax_CPAM54},
\cite{Lax_CPAM57}, \cite{OleiTMMO56} for the variational principle to one
quasilinear equation of the first order).

There were number of attempts to construct the generalized solutions to
system (\ref{1}) by analogy with 1-D case and with the aid of variational
principle, which can be written for 2-D inviscid Burgers equation and then
generalized through corresponding Hamilton-Jacobi equation. Further in
this paper we will show that in general this is impossible. The problems
for multi-D which involve different types of singularities can expose
rather complicated behavior, see, for example, \cite{MajdApMS84},
\cite{DiPeTAMS85}, \cite{ZhZhSIAM90}. Here we investigate the 'simplest'
degenerate case of 2-D system when the characteristics can be calculated
explicitly. Nevertheless the generalized solutions to system (\ref{1})
(see Definition 1 below) exhibit the different behavior compared to ones
for the single first order 2-D equation of Hamilton-Jacobi type. The main
ideas of the present paper were outlined in \cite{RykoKIAM30_98},
\cite{Ryko7HyC98}.

As it has been already said for smooth functions the system (\ref{1}) is
equivalent to (\ref{3}) and for small enough values of time and smooth
initial data the characteristics method can be applied. So the solution to
the problem (\ref{1}), (\ref{2}) reads
\begin{equation} \begin{aligned}
\ &\varrho(t,x,y)=\varrho_0(a,b)\frac{\partial (a,b)}{\partial (x,y)} \\
\ &u(t,x,y)=u(0,a,b)\quad ;\quad v(t,x,y)=v(0,a,b)\ ,
\end{aligned} \label{4} \end{equation}
where functions $a(t,x,y),b(t,x,y)$ are determined by the
equations
\begin{equation}
x=a+tu(0,a,b)\quad ; \quad y=b+tv(0,a,b)\ . \label{5}
\end{equation}
But it is well known that the characteristics of the system (\ref{3}) can
intersect themselves for finite time even for infinitely smooth initial
data.  Hence one finds the formation of singularities: discontinuities for
the velocity $(u,v)$ and $\delta$-functions for the density $\varrho$
which correspond to concentration process in some points or along some
curves. So it is necessary to introduce the notion of generalized solution
to the problem (\ref{1}), (\ref{2}) which is natural to formulate in terms
of Radon measures because of the type of singularities.

%
%  DEFINITION 1
%
\begin{defini}
Suppose $\bigl( P_t(dx,dy),I_t(dx,dy), J_t(dx,dy)\bigr)$ are the families
of Radon measures defined on Borel subsets of $\R^2$, weakly continuous
with respect to $t$ and such that $P_t \geq 0$, and $I_t,J_t$ are
absolutely continuous with respect to $P_t$ for almost every fixed $t>0$.
Let us define the vector function $\bigl( u(t,x,y),v(t,x,y) \bigr)$ as
Radon-Nykodim derivatives
\begin{equation}
u(t,x,y)=\frac{dI_t}{dP_t} \quad ; \quad
v(t,x,y)=\frac{dJ_t}{dP_t}\ .
\label{6} \end{equation}
Then $\bigl( P_t,I_t,J_t) \bigr)$ will be called the generalized
solution of the problem (\ref{1}), (\ref{2}), iff: 1) for an arbitrary
functions $f,g,h \in C_0^1(\R^2)$ (the space of continuously differentiable
functions with compact support) and $0<t_1<t_2<+\infty$
\vskip 6pt
\begin{align}
\ &\int\!\!\int f(x,y)P_{t_2}(dx,dy)-\int\!\!\int f(x,y)
P_{t_1}(dx,dy)= \notag \\
\ &\int_{t_1}^{t_2}\biggl\{
\int\!\!\int\frac{\partial f}{\partial x}(x,y)
I_{\tau}(dx,dy)+ \label{7} \\
\ &\int\!\!\int\frac{\partial f}{\partial y}
(x,y) J_{\tau}(dx,dy) \biggr\} d \tau \notag \\
\notag \\
\ &\int\!\!\int g(x,y)I_{t_2}(dx,dy)-\int\!\!\int g(x,y)
I_{t_1}(dx,dy)= \notag \\
\ &\int_{t_1}^{t_2}\biggl\{
\int\!\!\int\frac{\partial g}{\partial x}(x,y)u(\tau,x,y)
I_{\tau}(dx,dy)+ \label{8} \\
\ &\int\!\!\int\frac{\partial g}{\partial y}(x,y)v(\tau,x,y)
I_{\tau}(dx,dy) \biggr\} d \tau \notag \\
\notag \\
\ &\int\!\!\int h(x,y)J_{t_2}(dx,dy)-\int\!\!\int h(x,y)
J_{t_1}(dx,dy)= \notag \\
\ &\int_{t_1}^{t_2}\biggl\{
\int\!\!\int\frac{\partial h}{\partial x}(x,y)u(\tau,x,y)
J_{\tau}(dx,dy)+ \label{9} \\
\ &\int\!\!\int\frac{\partial h}{\partial y}(x,y)v(\tau,x,y)
J_{\tau}(dx,dy) \biggr\} d \tau\ , \notag
\end{align}
where $\int\!\!\int$ stands for the integration with respect to whole
$\R^2$;
\vskip 6pt
2) in the weak sense as $t \rightarrow +0$
\begin{align*}
\ &P_t \rightarrow \varrho(0,a,b)dadb \quad ; \quad
I_t \rightarrow \varrho(0,a,b)u(0,a,b)dadb\quad ;
\\
\ &J_t \rightarrow \varrho(0,a,b)v(0,a,b)dadb\quad.
\end{align*}
\end{defini}

Concerning the system (\ref{3}) there is well-known Hopf-Cole's
representation \cite{HopfCPAM50} for the solution with potential initial
velocity. This representation can be derived with the aid of infinitesimal
viscosity method (i.e. the adding of the second order operators with small
parameter $\varepsilon$ to the right-hand side of the last two equations
of (\ref{3}) and then letting $\varepsilon$ tend to zero) and allows to
determine the location of singularities. Namely
\begin{equation}
(u,v)=\nabla_{(x,y)}\Psi (t,x,y)\ ,\label{10}
\end{equation}
where
$$
\Psi (t,x,y)=\min\limits_{a,b}\ \bigl\{S_0(a,b)+
\frac{(x-a)^2}{2t}+\frac{(y-b)^2}{2t}\bigr\}\quad ,
$$
$S_0$ is the potential of initial velocity.
The representation (\ref{10}) is valid in domains of smoothness of $\Psi$
and gives the location of the set where singularities arise. That is the
singular points of $\Psi$ will be such points $(t,x,y)$ that the
global minimum of the expression in the braces with respect to
$(a,b)$ is attained more than in one point. For 1-D case the representation
(\ref{10}) in case of constant initial density gives us the way to
construct the generalized solution to the 1-D analogue of the problem
(\ref{1}), (\ref{2}), see \cite{ERSiUMNa95}, \cite{ERSiCMPh96} (in these
papers one can also find the generalization of (\ref{10}) to nonconstant
initial density). Unfortunately in 2-D case the generalized solution in
the sense of Definition 1 to the problem (\ref{1}), (\ref{2}) for constant
initial density  can not be constructed via (\ref{10}).

Our study is organized as follows.
In \S 2 we formulate the conditions on the initial data (\ref{2}) which
allow to conjecture the development of shock and then assuming that such
shock exists and satisfies the problem (\ref{1}), (\ref{2}) in the sense
of Definition 1 give the formulas for the velocity along the shock. These
formulas are the generalization of Rankine-Hugoniot conditions for
hyperbolic conservation laws. In \S 3 it is shown that in the case of
potential initial velocity these formulas determine the singularity
surface which in general does not coincide with one determined by the
representation (\ref{10}). As it is shown in \S 4 one can obtain
certain 1-D quasilinear hyperbolic system of PDEs for the motion {\it
inside} the shock. But this system does not satisfy Friedrichs'
symmetrizability condition, so it creates difficulties to prove the global
existence theorem. An interesting heuristic procedure to guess the system
is deferred to Appendix. Finally in \S 5 one investigates some simplest
solutions for derived system in case when $\varrho_0$, $u_0$, $v_0$ are
piecewise constant but the data along the initial shock front are rather
arbitrary. Even with such trivial external field one encounters with
complicate behavior of the generalized solutions. It is shown that there
are nontrivial cases when the solution can be defined explicitly and other
cases when the corresponding Cauchy problem is reduced to equation of type
$P_{tt}=const\cdot P_x$.  Such equations are ill-posed in the spaces of
functions of finite smoothness.

Further the letter subscripts (except '$i$', '$j$' and '$t$') will denote
appropriate derivatives.

%
%  SECTION 2
%
\section{The propagation of the shock front.}

Let us denote the coordinates in $\R^2\times\{t=0\}$ as $(a,b)$. Suppose
there is the $C^1$ curve $\G=\left(A(l),B(l)\right)$ in the
plane $(a,b)$, $l$ is a parameter, such that, for definiteness, $A_l\ge
0$, $B_l\ge 0$ ; $A_l^2+B_l^2\ne 0$ . Suppose that $\G$ can also
be found as the solution of the equation $G(a,b)=0$, $G\in C^1(\R^2)$.
Consider the following initial data (\ref{2})
\begin{align}
\varrho_0(a,b) &= \varrho_-+\left(\varrho_+-\varrho_-\right)H(G(a,b))
\notag \\
u_0(a,b) &= u_-+\left(u_+-u_-\right)H(G(a,b))\label{11}\\
v_0(a,b) &= v_-+\left(v_+-v_-\right)H(G(a,b))\ ,\notag
\end{align}
where $H$ is the Heaviside function,i.e. $H(\theta)=0\ \mbox{for}\ \theta<0,
H(\theta)=1\ \mbox{for}\ \theta>0$; and the following conditions hold
\medskip

I) $u_-(a,b), v_-(a,b)\in C^1(\overline{\G^-})\ ;\
u_+(a,b), v_+(a,b)\in C^1(\overline{\G^+})\ ,$

\hspace{3mm} here $\G^-\equiv \left\{(a,b)\ :\ G(a,b)<0\right\}$,
$\G^+\equiv \left\{(a,b)\ :\ G(a,b)>0\right\}$.
\smallskip

II) $A_lv_--B_lu_- < 0$ and $A_lv_+-B_lu_+ > 0$ on $\G$.

%
%  REMARK 2.1
%
\begin{remark}
Let us note that the condition I) reads that for every point
$(a^*,b^*)=(A(l^*),B(l^*))\in \G$ there exists some domain
$Q\subset \R^2\times [0,T(l^*)]$, $T(l^*)>0$, $(a^*,b^*)\in Q\bigcap
\{t=0\}$ where the transformations ($j=+,-$)
$$
\left\{ \begin{aligned}
\ &x=a+tu_j(a,b)
\\
\ &y=b+tv_j(a,b)
\end{aligned}\right.
$$
are nondegenerate for every $0<t<T(l^*)$, the characteristic lines
(\ref{5}) with $u(0,a,b)=u_-$, $v(0,a,b)=v_-$ never cross themselves, and
the characteristic lines (\ref{5}) with $u(0,a,b)=u_+$, $v(0,a,b)=v_+$
never cross themselves.

But the condition II) reads that for $0<t<T(l^*)$ there will be another
domain $Q_1\subset Q$ where exactly two characteristic lines issued from
different sides of $\G$ will come to the same point (among these
points the shock surface will form itself). So these conditions are
natural for the propagating shock front to exist locally up to time moment
$T(l^*)>0$.
\end{remark}

The following definition is essential for our analysis to elicit shock
fronts with 'good' behavior. Let us denote through ${\widehat Q}_1$ the
maximal $Q_1$-domain from Remark~2.1 which is valid for $l=l^*$.

%
%  DEFINITION 2
%
\begin{defini}
The shock front $\Sc$ will be called stable iff
$$\Sc\bigcap\{t=\tau\}\subset {\widehat Q}_1\bigcap\{t=\tau\} \quad
\mbox{as}\quad 0<\tau<T(l^*)\ .$$
\end{defini}

It is natural to seek the solution to the problem (\ref{1}), (\ref{11}) in
the form (see Fig.~\ref{fig1} below, {\it Eulerian} representation)
\begin{align}
\varrho(t,x,y) &= \varrho_- +
\left(\varrho_+-\varrho_-\right)H(S)+\widetilde {P_t}\delta (S)
\notag \\
\varrho u(t,x,y) &= \varrho_- u_- +
\left(\varrho_+ u_+-\varrho_- u_-\right)H(S)+\widetilde {I_t}
\delta (S)\label{12} \\
\varrho v(t,x,y) &= \varrho_- v_- +
\left(\varrho_+ v_+-\varrho_- v_-\right)H(S)+\widetilde {J_t}
\delta (S)\ ,\notag
\end{align}
where $\varrho_j(t,x,y),u_j(t,x,y),v_j(t,x,y),S(t,x,y)\in
C^1(\R_+\times \R^2)$, $(j=+,-)$ and satisfy (\ref{1}) in classical sense;
$S(t,x,y)=0$ represents some surface $\Sc$ in $\R_+\times \R^2$ and
$\Sc\bigcap \left(\{0\}\times \R^2\right)=\G$;
${\widetilde P}_t,{\widetilde I}_t,{\widetilde J}_t\in C^1(\Sc)$;
$H$ is the Heaviside function mentioned above, $\delta$ is standard Dirac
measure.

If the solution to the problem (\ref{1}), (\ref{11}) exists in the sense
of Definition 1 then by formulas (\ref{6}) there is defined the velocity
vector $\widetilde U=\left(\widetilde u,\widetilde v\right)$ on the
surface $\Sc$ and from each point of $\G$ one can draw the
integral curve of $\widetilde U$ on $\Sc$. Suppose
$\left(x^s(t,l),y^s(t,l)\right)$ is the corresponding parametrization of
$\Sc$ ($l$ is a parameter along $\G$) and
$d\left(x^s\right)/dt=\widetilde u$, $d\left(y^s\right)/dt=\widetilde v$.
Then according to formulas (\ref{5}) and conditions I), II) one can define
the functions $a_j(t,l), b_j(t,l)\in C^1(\Sc)$, $(j=+,-)$ which
correspond to the initial positions of characteristics which come to
the same points of the surface $\Sc$ from left and right.

Further let us define the map (see Fig.~\ref{fig1} below)
\begin{equation}
\Lc_t\ :\ (a,b)\rightarrow (x,y) \label{13}
\end{equation}
in the following way. Let us issue the characteristic line (keeping in mind
(\ref{11})) from the point $(a,b)$
\begin{equation}
x(\tau)=a+\tau u_0(a,b)\quad ;\quad y(\tau)=b+\tau v_0(a,b)
\label{14}
\end{equation}
and consider the time $\tau_0$ of the intersection of line (\ref{14}) and
$\Sc$. Then let us take $\Lc_t(a,b)=(x(t),y(t))$ if $\tau_0>t$
and $\Lc_t(a,b)=(x^s(t,l_0),y^s(t,l_0))$ if $\tau_0\le t$,
where $l_0$ is defined from the condition
$(x(\tau_0),y(\tau_0))=(x^s(\tau_0,l_0),y^s(\tau_0,l_0))$.

Now we are able to formulate the generalization of Rankine-Hugoniot
conditions.

%
%  THEOREM 2.1
%
\begin{theorem}
The generalized solution in the sense of Definition 1
to the problem (\ref{1}), (\ref{11}) in the form (\ref{12}) exists iff
the following formulas are true
\begin{equation}
\widetilde u=\widetilde {I_t}/\widetilde {P_t}\quad ;\quad
\widetilde v=\widetilde {J_t}/\widetilde {P_t}
\label{141}
\end{equation}
\begin{align}
\ \!\! &\widetilde {P_t}=\int\limits_0^t
\left[\varrho_0^+\left((a_+)_\tau (b_+)_l-(b_+)_\tau (a_+)_l\right)-
\varrho_0^-\left((a_-)_\tau (b_-)_l-(b_-)_\tau (a_-)_l\right)\right]d\tau
\notag \\
\ \!\! &\widetilde {I_t}=\int\limits_0^t
\left[\varrho_0^+u_0^+\left((a_+)_\tau (b_+)_l-(b_+)_\tau (a_+)_l\right)-
\varrho_0^-u_0^-\left((a_-)_\tau (b_-)_l-(b_-)_\tau (a_-)_l\right)\right]
d\tau
\label{15} \\
\ \!\! &\widetilde {J_t}=\int\limits_0^t
\left[\varrho_0^+v_0^+\left((a_+)_\tau (b_+)_l-(b_+)_\tau (a_+)_l\right)-
\varrho_0^-v_0^-\left((a_-)_\tau (b_-)_l-(b_-)_\tau (a_-)_l\right)\right]
d\tau , \notag
\end{align}
where $\varrho_0^j\equiv\varrho_0(a_j,b_j)$, $u_0^j\equiv u_0(a_j,b_j)$,
$v_0^j\equiv v_0(a_j,b_j)$, $(j=+,-)$.
\end{theorem}

{\sl PROOF.} Suppose there are some family of domains $D(\tau)$ with
oriented boundaries $\partial D(\tau)$ in the plane $(a,b)$,
$D(\tau_1)\subset D(\tau_2)$ for $0<\tau_1<\tau_2<T$.  Suppose $\partial
D(\tau)$ is closed curve $(a(\tau,l), b(\tau,l))$, $l$ is a parameter
along the curve, $a,b\in C^1([\tau_1,\tau_2] \times \R)$.  Then it is easy
to check that the following formula is valid
\begin{equation}
\frac{d}{d\tau}\int\!\!\int\limits_{D(\tau)}\varphi dadb=
\int\!\!\int\limits_{D(\tau)}\frac{\partial\varphi}{\partial\tau}dadb+
\oint\limits_{\partial D(\tau)}\varphi (a_\tau b_l-b_\tau a_l)dl\quad ,
\label{16}
\end{equation}
where $\varphi\in C^1\left([\tau_1,\tau_2]\times\overline{\bigcup_{\tau\in
[\tau_1,\tau_2]}D(\tau)}\right)$.

Denote through $\Pi$ some circle in the plane $(a,b)$ such that
$\Pi\times\{t\}\bigcap\Sc_t\ne\emptyset$,
$\Sc_t\equiv\Sc\bigcap\{t\}\times \R^2$.
Denote through $D_-(t)$, $D_+(t)$ (Consult with the Fig.~\ref{fig1}.) the
domains in the plain $(x,y)$ to which at time $t$ the characteristics come
from $(a_-,b_-)$, $(a_+,b_+)$ respectively (note that $D_-(t)\bigcup
D_+(t)=\Pi\times\{t\}$ and $D_-(t)\bigcap D_+(t)=\Sc_t$).  Further let us
consider the right hand side of integral identity (\ref{7}).  Take $f\in
C_0^1(\R^2)$ such that $supp\ f\subset\Pi$. Suppose $D_j^*(t)=
\Lc_t^{-1}(D_j(t))$, $(j=+,-)$, $D_\Sc(t)=\Lc_t^{-1}(\Sc_t)$.  Define the
function $f_t^*(a,b)$ such that $f_t^*(a,b)= f(\Lc_t(a,b))$. Then one has
\begin{figure}[tph]
%
%TexCad Options
%\grade{\on}
%\emlines{\on}
%\beziermacro{\off}
%\reduce{\on}
%\snapping{\off}
%\quality{2.00}
%\graddiff{0.01}
%\snapasp{1}
%\zoom{1.00}
\special{em:linewidth 0.4pt}
\unitlength .4mm
\linethickness{0.4pt}
\hspace{3cm}
\begin{picture}(124.33,152.67)
%\bezier{316}(16.33,81.67)(16.67,123.00)(47.67,101.67)
\emline{16.33}{81.67}{1}{16.41}{84.87}{2}
\emline{16.41}{84.87}{3}{16.58}{87.86}{4}
\emline{16.58}{87.86}{5}{16.85}{90.66}{6}
\emline{16.85}{90.66}{7}{17.22}{93.26}{8}
\emline{17.22}{93.26}{9}{17.69}{95.67}{10}
\emline{17.69}{95.67}{11}{18.25}{97.87}{12}
\emline{18.25}{97.87}{13}{18.91}{99.87}{14}
\emline{18.91}{99.87}{15}{19.67}{101.67}{16}
\emline{19.67}{101.67}{17}{20.53}{103.27}{18}
\emline{20.53}{103.27}{19}{21.48}{104.68}{20}
\emline{21.48}{104.68}{21}{22.53}{105.88}{22}
\emline{22.53}{105.88}{23}{23.68}{106.88}{24}
\emline{23.68}{106.88}{25}{24.93}{107.69}{26}
\emline{24.93}{107.69}{27}{26.28}{108.29}{28}
\emline{26.28}{108.29}{29}{27.72}{108.70}{30}
\emline{27.72}{108.70}{31}{29.26}{108.90}{32}
\emline{29.26}{108.90}{33}{30.90}{108.91}{34}
\emline{30.90}{108.91}{35}{32.63}{108.71}{36}
\emline{32.63}{108.71}{37}{34.47}{108.32}{38}
\emline{34.47}{108.32}{39}{36.40}{107.73}{40}
\emline{36.40}{107.73}{41}{38.43}{106.93}{42}
\emline{38.43}{106.93}{43}{40.55}{105.94}{44}
\emline{40.55}{105.94}{45}{42.78}{104.75}{46}
\emline{42.78}{104.75}{47}{45.10}{103.36}{48}
\emline{45.10}{103.36}{49}{47.67}{101.67}{50}
%\end
%\bezier{232}(47.67,101.67)(69.67,78.00)(49.00,62.00)
\emline{47.67}{101.67}{51}{49.63}{99.48}{52}
\emline{49.63}{99.48}{53}{51.40}{97.32}{54}
\emline{51.40}{97.32}{55}{52.99}{95.20}{56}
\emline{52.99}{95.20}{57}{54.39}{93.11}{58}
\emline{54.39}{93.11}{59}{55.61}{91.05}{60}
\emline{55.61}{91.05}{61}{56.64}{89.03}{62}
\emline{56.64}{89.03}{63}{57.48}{87.04}{64}
\emline{57.48}{87.04}{65}{58.14}{85.08}{66}
\emline{58.14}{85.08}{67}{58.62}{83.15}{68}
\emline{58.62}{83.15}{69}{58.91}{81.26}{70}
\emline{58.91}{81.26}{71}{59.01}{79.41}{72}
\emline{59.01}{79.41}{73}{58.93}{77.58}{74}
\emline{58.93}{77.58}{75}{58.66}{75.79}{76}
\emline{58.66}{75.79}{77}{58.21}{74.04}{78}
\emline{58.21}{74.04}{79}{57.57}{72.31}{80}
\emline{57.57}{72.31}{81}{56.75}{70.62}{82}
\emline{56.75}{70.62}{83}{55.74}{68.97}{84}
\emline{55.74}{68.97}{85}{54.54}{67.34}{86}
\emline{54.54}{67.34}{87}{53.16}{65.75}{88}
\emline{53.16}{65.75}{89}{51.60}{64.20}{90}
\emline{51.60}{64.20}{91}{49.00}{62.00}{92}
%\end
%\bezier{264}(49.00,62.00)(16.00,50.33)(16.33,81.00)
\emline{49.00}{62.00}{93}{46.18}{61.06}{94}
\emline{46.18}{61.06}{95}{43.49}{60.28}{96}
\emline{43.49}{60.28}{97}{40.92}{59.67}{98}
\emline{40.92}{59.67}{99}{38.49}{59.22}{100}
\emline{38.49}{59.22}{101}{36.17}{58.92}{102}
\emline{36.17}{58.92}{103}{33.99}{58.79}{104}
\emline{33.99}{58.79}{105}{31.93}{58.82}{106}
\emline{31.93}{58.82}{107}{30.01}{59.02}{108}
\emline{30.01}{59.02}{109}{28.20}{59.37}{110}
\emline{28.20}{59.37}{111}{26.53}{59.89}{112}
\emline{26.53}{59.89}{113}{24.98}{60.56}{114}
\emline{24.98}{60.56}{115}{23.56}{61.40}{116}
\emline{23.56}{61.40}{117}{22.27}{62.40}{118}
\emline{22.27}{62.40}{119}{21.10}{63.56}{120}
\emline{21.10}{63.56}{121}{20.07}{64.89}{122}
\emline{20.07}{64.89}{123}{19.15}{66.37}{124}
\emline{19.15}{66.37}{125}{18.37}{68.02}{126}
\emline{18.37}{68.02}{127}{17.71}{69.83}{128}
\emline{17.71}{69.83}{129}{17.18}{71.80}{130}
\emline{17.18}{71.80}{131}{16.78}{73.93}{132}
\emline{16.78}{73.93}{133}{16.51}{76.22}{134}
\emline{16.51}{76.22}{135}{16.36}{78.68}{136}
\emline{16.36}{78.68}{137}{16.33}{81.00}{138}
%\end
%\bezier{540}(72.00,83.00)(83.33,152.67)(104.00,91.67)
\emline{72.00}{83.00}{139}{72.70}{87.13}{140}
\emline{72.70}{87.13}{141}{73.42}{91.01}{142}
\emline{73.42}{91.01}{143}{74.15}{94.65}{144}
\emline{74.15}{94.65}{145}{74.90}{98.05}{146}
\emline{74.90}{98.05}{147}{75.67}{101.20}{148}
\emline{75.67}{101.20}{149}{76.46}{104.11}{150}
\emline{76.46}{104.11}{151}{77.26}{106.78}{152}
\emline{77.26}{106.78}{153}{78.08}{109.21}{154}
\emline{78.08}{109.21}{155}{78.92}{111.39}{156}
\emline{78.92}{111.39}{157}{79.78}{113.33}{158}
\emline{79.78}{113.33}{159}{80.65}{115.03}{160}
\emline{80.65}{115.03}{161}{81.54}{116.48}{162}
\emline{81.54}{116.48}{163}{82.45}{117.70}{164}
\emline{82.45}{117.70}{165}{83.37}{118.67}{166}
\emline{83.37}{118.67}{167}{84.32}{119.39}{168}
\emline{84.32}{119.39}{169}{85.28}{119.87}{170}
\emline{85.28}{119.87}{171}{86.25}{120.11}{172}
\emline{86.25}{120.11}{173}{87.25}{120.11}{174}
\emline{87.25}{120.11}{175}{88.26}{119.87}{176}
\emline{88.26}{119.87}{177}{89.29}{119.38}{178}
\emline{89.29}{119.38}{179}{90.34}{118.65}{180}
\emline{90.34}{118.65}{181}{91.40}{117.67}{182}
\emline{91.40}{117.67}{183}{92.48}{116.46}{184}
\emline{92.48}{116.46}{185}{93.58}{115.00}{186}
\emline{93.58}{115.00}{187}{94.70}{113.29}{188}
\emline{94.70}{113.29}{189}{95.83}{111.35}{190}
\emline{95.83}{111.35}{191}{96.98}{109.16}{192}
\emline{96.98}{109.16}{193}{98.15}{106.73}{194}
\emline{98.15}{106.73}{195}{99.33}{104.05}{196}
\emline{99.33}{104.05}{197}{100.54}{101.14}{198}
\emline{100.54}{101.14}{199}{101.76}{97.98}{200}
\emline{101.76}{97.98}{201}{104.00}{91.67}{202}
%\end
%\bezier{296}(104.00,91.67)(114.67,44.33)(90.67,53.00)
\emline{104.00}{91.67}{203}{104.82}{87.86}{204}
\emline{104.82}{87.86}{205}{105.52}{84.24}{206}
\emline{105.52}{84.24}{207}{106.11}{80.81}{208}
\emline{106.11}{80.81}{209}{106.58}{77.57}{210}
\emline{106.58}{77.57}{211}{106.92}{74.53}{212}
\emline{106.92}{74.53}{213}{107.16}{71.67}{214}
\emline{107.16}{71.67}{215}{107.27}{69.01}{216}
\emline{107.27}{69.01}{217}{107.26}{66.53}{218}
\emline{107.26}{66.53}{219}{107.14}{64.25}{220}
\emline{107.14}{64.25}{221}{106.90}{62.15}{222}
\emline{106.90}{62.15}{223}{106.54}{60.25}{224}
\emline{106.54}{60.25}{225}{106.07}{58.53}{226}
\emline{106.07}{58.53}{227}{105.47}{57.01}{228}
\emline{105.47}{57.01}{229}{104.76}{55.68}{230}
\emline{104.76}{55.68}{231}{103.93}{54.54}{232}
\emline{103.93}{54.54}{233}{102.98}{53.58}{234}
\emline{102.98}{53.58}{235}{101.92}{52.82}{236}
\emline{101.92}{52.82}{237}{100.73}{52.25}{238}
\emline{100.73}{52.25}{239}{99.43}{51.87}{240}
\emline{99.43}{51.87}{241}{98.01}{51.68}{242}
\emline{98.01}{51.68}{243}{96.48}{51.68}{244}
\emline{96.48}{51.68}{245}{94.82}{51.87}{246}
\emline{94.82}{51.87}{247}{93.05}{52.26}{248}
\emline{93.05}{52.26}{249}{90.67}{53.00}{250}
%\end
%\bezier{172}(90.67,53.00)(68.67,66.00)(72.00,82.67)
\emline{90.67}{53.00}{251}{88.19}{54.52}{252}
\emline{88.19}{54.52}{253}{85.89}{56.07}{254}
\emline{85.89}{56.07}{255}{83.76}{57.65}{256}
\emline{83.76}{57.65}{257}{81.80}{59.24}{258}
\emline{81.80}{59.24}{259}{80.02}{60.87}{260}
\emline{80.02}{60.87}{261}{78.40}{62.52}{262}
\emline{78.40}{62.52}{263}{76.96}{64.19}{264}
\emline{76.96}{64.19}{265}{75.68}{65.89}{266}
\emline{75.68}{65.89}{267}{74.58}{67.61}{268}
\emline{74.58}{67.61}{269}{73.65}{69.36}{270}
\emline{73.65}{69.36}{271}{72.89}{71.13}{272}
\emline{72.89}{71.13}{273}{72.30}{72.92}{274}
\emline{72.30}{72.92}{275}{71.88}{74.75}{276}
\emline{71.88}{74.75}{277}{71.64}{76.59}{278}
\emline{71.64}{76.59}{279}{71.56}{78.46}{280}
\emline{71.56}{78.46}{281}{71.66}{80.36}{282}
\emline{71.66}{80.36}{283}{72.00}{82.67}{284}
%\end
%\bezier{208}(16.33,81.67)(44.33,95.00)(59.00,80.33)
\emline{16.33}{81.67}{285}{19.06}{82.91}{286}
\emline{19.06}{82.91}{287}{21.72}{84.02}{288}
\emline{21.72}{84.02}{289}{24.32}{85.00}{290}
\emline{24.32}{85.00}{291}{26.85}{85.83}{292}
\emline{26.85}{85.83}{293}{29.32}{86.54}{294}
\emline{29.32}{86.54}{295}{31.72}{87.10}{296}
\emline{31.72}{87.10}{297}{34.06}{87.53}{298}
\emline{34.06}{87.53}{299}{36.33}{87.83}{300}
\emline{36.33}{87.83}{301}{38.54}{87.99}{302}
\emline{38.54}{87.99}{303}{40.69}{88.01}{304}
\emline{40.69}{88.01}{305}{42.76}{87.89}{306}
\emline{42.76}{87.89}{307}{44.78}{87.65}{308}
\emline{44.78}{87.65}{309}{46.73}{87.26}{310}
\emline{46.73}{87.26}{311}{48.61}{86.74}{312}
\emline{48.61}{86.74}{313}{50.43}{86.08}{314}
\emline{50.43}{86.08}{315}{52.19}{85.29}{316}
\emline{52.19}{85.29}{317}{53.88}{84.36}{318}
\emline{53.88}{84.36}{319}{55.51}{83.30}{320}
\emline{55.51}{83.30}{321}{57.07}{82.10}{322}
\emline{57.07}{82.10}{323}{59.00}{80.33}{324}
%\end
%\bezier{144}(73.67,91.67)(89.33,102.00)(103.67,92.33)
\emline{73.67}{91.67}{325}{75.84}{93.01}{326}
\emline{75.84}{93.01}{327}{77.99}{94.15}{328}
\emline{77.99}{94.15}{329}{80.14}{95.10}{330}
\emline{80.14}{95.10}{331}{82.27}{95.86}{332}
\emline{82.27}{95.86}{333}{84.39}{96.43}{334}
\emline{84.39}{96.43}{335}{86.49}{96.81}{336}
\emline{86.49}{96.81}{337}{88.58}{96.99}{338}
\emline{88.58}{96.99}{339}{90.66}{96.98}{340}
\emline{90.66}{96.98}{341}{92.73}{96.77}{342}
\emline{92.73}{96.77}{343}{94.78}{96.37}{344}
\emline{94.78}{96.37}{345}{96.82}{95.78}{346}
\emline{96.82}{95.78}{347}{98.85}{95.00}{348}
\emline{98.85}{95.00}{349}{100.87}{94.02}{350}
\emline{100.87}{94.02}{351}{103.67}{92.33}{352}
%\end
%\bezier{176}(73.33,70.00)(100.33,82.33)(107.33,69.33)
\emline{73.33}{70.00}{353}{76.34}{71.32}{354}
\emline{76.34}{71.32}{355}{79.21}{72.48}{356}
\emline{79.21}{72.48}{357}{81.96}{73.47}{358}
\emline{81.96}{73.47}{359}{84.57}{74.30}{360}
\emline{84.57}{74.30}{361}{87.06}{74.96}{362}
\emline{87.06}{74.96}{363}{89.42}{75.46}{364}
\emline{89.42}{75.46}{365}{91.65}{75.80}{366}
\emline{91.65}{75.80}{367}{93.75}{75.98}{368}
\emline{93.75}{75.98}{369}{95.72}{75.99}{370}
\emline{95.72}{75.99}{371}{97.56}{75.84}{372}
\emline{97.56}{75.84}{373}{99.27}{75.52}{374}
\emline{99.27}{75.52}{375}{100.85}{75.04}{376}
\emline{100.85}{75.04}{377}{102.31}{74.40}{378}
\emline{102.31}{74.40}{379}{103.63}{73.59}{380}
\emline{103.63}{73.59}{381}{104.83}{72.62}{382}
\emline{104.83}{72.62}{383}{105.90}{71.49}{384}
\emline{105.90}{71.49}{385}{107.33}{69.33}{386}
%\end
%\vector(43.67,40.00)(80.67,40.00)
\put(80.67,40.00){\vector(1,0){0.2}}
\emline{43.67}{40.00}{387}{80.67}{40.00}{388}
%\end
%\bezvec{148}(123.33,91.67)(114.00,78.67)(98.33,92.33)
\put(98.33,92.33){\vector(-4,3){0.2}}
\emline{123.33}{91.67}{389}{122.04}{90.03}{390}
\emline{122.04}{90.03}{391}{120.70}{88.64}{392}
\emline{120.70}{88.64}{393}{119.29}{87.49}{394}
\emline{119.29}{87.49}{395}{117.83}{86.59}{396}
\emline{117.83}{86.59}{397}{116.30}{85.93}{398}
\emline{116.30}{85.93}{399}{114.72}{85.51}{400}
\emline{114.72}{85.51}{401}{113.09}{85.33}{402}
\emline{113.09}{85.33}{403}{111.39}{85.40}{404}
\emline{111.39}{85.40}{405}{109.64}{85.72}{406}
\emline{109.64}{85.72}{407}{107.83}{86.27}{408}
\emline{107.83}{86.27}{409}{105.96}{87.07}{410}
\emline{105.96}{87.07}{411}{104.03}{88.12}{412}
\emline{104.03}{88.12}{413}{102.05}{89.40}{414}
\emline{102.05}{89.40}{415}{98.33}{92.33}{416}
%\end
%\bezvec{112}(124.33,72.00)(121.67,81.33)(104.33,76.00)
\put(104.33,76.00){\vector(-4,-1){0.2}}
\emline{124.33}{72.00}{417}{123.74}{73.55}{418}
\emline{123.74}{73.55}{419}{122.91}{74.87}{420}
\emline{122.91}{74.87}{421}{121.85}{75.95}{422}
\emline{121.85}{75.95}{423}{120.56}{76.80}{424}
\emline{120.56}{76.80}{425}{119.03}{77.41}{426}
\emline{119.03}{77.41}{427}{117.27}{77.79}{428}
\emline{117.27}{77.79}{429}{115.27}{77.94}{430}
\emline{115.27}{77.94}{431}{113.04}{77.85}{432}
\emline{113.04}{77.85}{433}{110.58}{77.53}{434}
\emline{110.58}{77.53}{435}{107.88}{76.97}{436}
\emline{107.88}{76.97}{437}{104.33}{76.00}{438}
%\end
\put(155.00,98.67){\makebox(0,0)[cc]{$\left(a_-(t,l),b_-(t,l)\right)$}}
\put(155.33,65.00){\makebox(0,0)[cc]{$\left(a_+(t,l),b_+(t,l)\right)$}}
\put(61.33,45.00){\makebox(0,0)[cc]{$l$}}
\put(88.00,131.33){\makebox(0,0)[cc]{{\it Lagrangian}}}
\put(35.33,130.67){\makebox(0,0)[cc]{{\it Eulerian}}}
\put(35.33,95.67){\makebox(0,0)[cc]{$D_-$}}
\put(32.67,71.33){\makebox(0,0)[cc]{$D_+$}}
\put(51.33,78.33){\makebox(0,0)[cc]{$\bf\Sc$}}
\put(87.33,105.33){\makebox(0,0)[cc]{$D_-^*$}}
\put(88.00,83.67){\makebox(0,0)[cc]{$D_{\bf\Sc}$}}
\put(93.00,65.00){\makebox(0,0)[cc]{$D_+^*$}}
\end{picture}
\vspace{-1.5cm}
\caption{\it The map $\Lc_t$ for fixed $t$; $l$ is a parameter along
$\Sc$.}
\label{fig1}
\end{figure}

\begin{align*}
\ &\int\limits_{t_1}^{t_2}dt\left\{\sum\limits_{j=+,-}
\int\!\!\int\limits_{D_j(t)}\left(f_x u+f_y v\right)P_t(dx,dy)
+\int\limits_{{\cal S}_t}\left(f_x u+f_y v\right)P_t(dl)\right\}= \\
\ &\int\limits_{t_1}^{t_2}dt\Biggl\{
\sum\limits_{j=+,-}\int\!\!\int\limits_{D_j^*(t)}
\Bigl[f_x\left(a+t u_0,b+t v_0\right)\varrho_0(a,b)u_0(a,b)+ \\
\ &f_y\left(a+t u_0,b+t v_0\right)\varrho_0(a,b)v_0(a,b)\Bigr]dadb+
\int\limits_{\Sc_t}
\Bigl[f_x\left(x^s(t,l),y^s(t,l)\right)u(t,l)+ \\
\ &f_y\left(x^s(t,l),y^s(t,l)\right)v(t,l)\Bigr]P_t(dl)\Biggr\}=
\int\limits_{t_1}^{t_2}dt\Biggl\{
\sum\limits_{j=+,-}\int\!\!\int\limits_{D_j^*(t)}
\frac{\partial f_t^*}{\partial t}\varrho_0(a,b)dadb+ \\
\ &\int\limits_{\Sc_t}\frac{\partial f_t^*}{\partial t}
\left(x^s(t,l),y^s(t,l)\right)P_t(dl)\Biggr\}=\int\limits_{t_1}^{t_2}dt
\Biggl\{\sum\limits_{j=+,-}\frac{d}{dt}\int\!\!\int\limits_{D_j^*(t)}f_t^*
\varrho_0(a,b)dadb- \\
\ &\int\limits_{\Sc_t}
f\left(x^s(t,l),y^s(t,l)\right)\biggl[\varrho_0(a_-,b_-)
\left(\frac{\partial (a_-)}{\partial t}\frac{\partial (b_-)}{\partial l}-
\frac{\partial (b_-)}{\partial t}\frac{\partial (a_-)}{\partial l}\right)-
\\
\ &\varrho_0(a_+,b_+)
\left(\frac{\partial (a_+)}{\partial t}\frac{\partial (b_+)}{\partial l}-
\frac{\partial (b_+)}{\partial t}\frac{\partial (a_+)}{\partial l}\right)
\biggr]dl+ \\
\ &\int\limits_{\Sc_t}\frac{df_t^*}{dt}\left(x^s(t,l),y^s(t,l)\right)
P_t(dl)\Biggr\}\ ,
\end{align*}
using formula (\ref{16}).

Consider the right hand side of integral identity (\ref{8}), here $g$ is
playing the role of $f$,
\begin{align*}
\ &\int\limits_{t_1}^{t_2}dt\left\{\sum\limits_{j=+,-}
\int\!\!\int\limits_{D_j(t)}\left(g_x u+g_y v\right)I_t(dx,dy)
+\int\limits_{{\cal S}_t}\left(g_x u+g_y v\right)I_t(dl)\right\}= \\
\ &\int\limits_{t_1}^{t_2}dt\Biggl\{
\sum\limits_{j=+,-}\int\!\!\int\limits_{D_j^*(t)}
\Bigl[g_x\left(a+t u_0,b+t v_0\right)\varrho_0(a,b)u_0^2(a,b)+ \\
\ &
g_y\left(a+t u_0,b+t v_0\right)\varrho_0(a,b)u_0(a,b) v_0(a,b)\Bigr]dadb+
\\
\ &\int\limits_{\Sc_t}
\Bigl[g_x\left(x^s(t,l),y^s(t,l)\right)u(t,l)+
g_y\left(x^s(t,l),y^s(t,l)\right)v(t,l)\Bigr]I_t(dl)\Biggr\}= \\
\ &\int\limits_{t_1}^{t_2}dt\Biggl\{
\sum\limits_{j=+,-}\int\!\!\int\limits_{D_j^*(t)}
\frac{\partial g_t^*}{\partial t}\varrho_0(a,b)u_0(a,b)dadb+
\int\limits_{\Sc_t}\frac{\partial g_t^*}{\partial t}
\left(x^s(t,l),y^s(t,l)\right) I_t(dl)\Biggr\}= \\
\ &\int\limits_{t_1}^{t_2}dt\Biggl\{
\sum\limits_{j=+,-}\frac{d}{dt}\int\!\!\int\limits_{D_j^*(t)}g_t^*
\varrho_0(a,b)u_0(a,b)dadb- \\
\ &\int\limits_{\Sc_t}
g\left(x^s(t,l),y^s(t,l)\right)\biggl[\varrho_0^-u_0^-
\left(\frac{\partial (a_-)}{\partial t}\frac{\partial (b_-)}{\partial l}-
\frac{\partial (b_-)}{\partial t}\frac{\partial (a_-)}{\partial l}\right)-
\\
\ &\varrho_0^+u_0^+
\left(\frac{\partial (a_+)}{\partial t}\frac{\partial (b_+)}{\partial l}-
\frac{\partial (b_+)}{\partial t}\frac{\partial (a_+)}{\partial l}\right)
\biggr]dl+
\int\limits_{\Sc_t}\frac{\partial g_t^*}{\partial t}
\left(x^s(t,l),y^s(t,l)\right)I_t(dl)\Biggr\},
\end{align*}
using formula (\ref{16}).

Consider the right hand side of integral identity (\ref{9}), here $h$ is
playing the role of $f$,
\begin{align*}
\ &\int\limits_{t_1}^{t_2}dt\left\{\sum\limits_{j=+,-}
\int\!\!\int\limits_{D_j(t)}\left(h_x u+h_y v\right)J_t(dx,dy)
+\int\limits_{{\cal S}_t}\left(h_x u+h_y v\right)J_t(dl)\right\}= \\
\ &\int\limits_{t_1}^{t_2}dt\Biggl\{
\sum\limits_{j=+,-}\int\!\!\int\limits_{D_j^*(t)}
\Bigl[h_x\left(a+t u_0,b+t v_0\right)\varrho_0(a,b)u_0(a,b) v_0(a,b)+ \\
\ &h_y\left(a+t u_0,b+t v_0\right)\varrho_0(a,b)v_0^2(a,b)\Bigr]dadb+ \\
\ &\int\limits_{\Sc_t}
\Bigl[h_x\left(x^s(t,l),y^s(t,l)\right)u(t,l)+
h_y\left(x^s(t,l),y^s(t,l)\right)v(t,l)\Bigr]J_t(dl)\Biggr\}= \\
\ &\int\limits_{t_1}^{t_2}dt\Biggl\{
\sum\limits_{j=+,-}\int\!\!\int\limits_{D_j^*(t)}
\frac{\partial h_t^*}{\partial t}\varrho_0(a,b)v_0(a,b)dadb+
\int\limits_{\Sc_t}\frac{\partial h_t^*}{\partial t}
\left(x^s(t,l),y^s(t,l)\right)
J_t(dl)\Biggr\}= \\
\ &\int\limits_{t_1}^{t_2}dt\Biggl\{
\sum\limits_{j=+,-}\frac{d}{dt}\int\!\!\int\limits_{D_j^*(t)}h_t^*
\varrho_0(a,b)v_0(a,b)dadb- \\
\ &\int\limits_{\Sc_t}
h\left(x^s(t,l),y^s(t,l)\right)\biggl[\varrho_0^-v_0^-
\left(\frac{\partial (a_-)}{\partial t}\frac{\partial (b_-)}{\partial l}-
\frac{\partial (b_-)}{\partial t}\frac{\partial (a_-)}{\partial l}\right)-
\\
\ &\varrho_0^+v_0^+
\left(\frac{\partial (a_+)}{\partial t}\frac{\partial (b_+)}{\partial l}-
\frac{\partial (b_+)}{\partial t}\frac{\partial (a_+)}{\partial l}\right)
\biggr]dl+
\int\limits_{\Sc_t}\frac{\partial h_t^*}{\partial t}
\left(x^s(t,l),y^s(t,l)\right)J_t(dl)\Biggr\},
\end{align*}
again using formula (\ref{16}).

\epr

%
%  THEOREM 2.2
%
\begin{theorem}
For every $l$ and $0<t<T(l)$ the following formulas are true
\begin{align}
\ &\int\limits_{0}^{t}\left[x^s(t,l)-a_+(\tau,l)-tu_0^+(\tau,l)\right]
\varrho_0(a_+,b_+)\left((a_+)_\tau (b_+)_l-(b_+)_\tau (a_+)_l\right)d\tau=
\notag \\
\ &\int\limits_{0}^{t}\left[x^s(t,l)-a_-(\tau,l)-tu_0^-(\tau,l)\right]
\varrho_0(a_-,b_-)\left((a_-)_\tau (b_-)_l-(b_-)_\tau (a_-)_l\right)d\tau
\ ; \label{17} \\
\ &\int\limits_{0}^{t}\left[y^s(t,l)-b_+(\tau,l)-tv_0^+(\tau,l)\right]
\varrho_0(a_+,b_+)\left((a_+)_\tau (b_+)_l-(b_+)_\tau (a_+)_l\right)d\tau=
\notag \\
\ &\int\limits_{0}^{t}\left[y^s(t,l)-b_-(\tau,l)-tv_0^-(\tau,l)\right]
\varrho_0(a_-,b_-)\left((a_-)_\tau (b_-)_l-(b_-)_\tau (a_-)_l\right)d\tau
\notag
\end{align}
\end{theorem}

{\sl PROOF.} Let us take formulas (\ref{141}) and write them in the
following way
$$
\left(x^s\right)_\tau P_{t=\tau}=I_{t=\tau}\quad ;\quad
\left(y^s\right)_\tau P_{t=\tau}=J_{t=\tau}\ ,
$$
Now let us integrate these equalities with respect to $\tau$ from $0$
to $t$ and then integrating by parts in both sides one obtains
\begin{align*}
\ &\int\limits_{0}^{t}\left(x^s(t)-x^s(\tau)\right)\times \\
\ &\quad
\left[\varrho_0^+\left((a_+)_\tau (b_+)_l-(b_+)_\tau (a_+)_l\right)-
\varrho_0^-\left((a_-)_\tau (b_-)_l-(b_-)_\tau (a_-)_l\right)\right]d\tau=
\\
\ &\int\limits_{0}^{t}\left(t-\tau\right)\times \\
\ &\quad
\left[\varrho_0^+u_0^+\left((a_+)_\tau (b_+)_l-(b_+)_\tau (a_+)_l\right)-
\varrho_0^-u_0^-\left((a_-)_\tau (b_-)_l-(b_-)_\tau (a_-)_l\right)
\right]d\tau
\\
\ &\int\limits_{0}^{t}\left(y^s(t)-y^s(\tau)\right)\times \\
\ &\quad
\left[\varrho_0^+\left((a_+)_\tau (b_+)_l-(b_+)_\tau (a_+)_l\right)-
\varrho_0^-\left((a_-)_\tau (b_-)_l-(b_-)_\tau (a_-)_l\right)\right]d\tau=
\\
\ &\int\limits_{0}^{t}\left(t-\tau\right)\times \\
\ &\quad
\left[\varrho_0^+v_0^+\left((a_+)_\tau (b_+)_l-(b_+)_\tau (a_+)_l\right)-
\varrho_0^-v_0^-\left((a_-)_\tau (b_-)_l-(b_-)_\tau (a_-)_l\right)
\right]d\tau .
\end{align*}
Then applying formulas (\ref{14}) we come to the assertion of Theorem 2.2.

\epr

%
%  REMARK 2.2
%
\begin{remark}
Let us mention that expressions (\ref{17}) look similar to the expressions
for the adhesion principle in \cite{ERSiUMNa95}, \cite{ERSiCMPh96}. So all
characteristics initially started at points $(a_j(\tau,l),b_j(\tau,l))$,
$j=+,-$, $\tau <t$, $l$ is fixed, will 'concentrate' in one point at time
$t$.
\end{remark}

%
%  SECTION 3
%
\section{Comparison with the variational representation.}

In this section we assume that $\varrho_0(a,b)\equiv 1$.
Now let us compare the formulas (\ref{15}) for shocks with the restrictions
on singularity surface which can be obtained through the variational
representation (\ref{10}) in case of potential initial data. Suppose
$\Fc$=$(x^v(t,l),y^v(t,l))$ is the singularity surface obtained from
(\ref{10}). Then $\Fc$ can be represented as the solution of the equation
\begin{equation}
F(a_+,b_+;t,x,y)-F(a_-,b_-;t,x,y)=0\quad ,\label{18}
\end{equation}
where
$$F\equiv S_0(a,b)+\frac{(x-a)^2}{2t}+\frac{(y-b)^2}{2t}\quad ,$$
and $a_j(t,x,y),b_j(t,x,y)$, $(j=+,-)$ are defined from the
system of equations
$$F_a(a,b;t,x,y)=0\quad ;\quad F_b(a,b;t,x,y)=0\quad .$$

%
%  THEOREM 3.1
%
\begin{theorem}
Suppose the initial velocity vector in (\ref{11}) is a potential vector and
the potential $S_0(a,b)$ has the form
\begin{equation}\left\{\begin{aligned}
\ &S_0(a,b)=0\ ,\ \mbox{as}\ b>-\varepsilon f(a,b) \\
\ &S_0(a,b)=b+\varepsilon f(a,b)\ ,\ \mbox{as}\ b<-\varepsilon f(a,b)\ ,
\end{aligned}\right.\label{19}
\end{equation}
where $\varepsilon>0$ is sufficiently small and the function
$f\in C^\infty (\R^2)$, $f(0,0)=0$, satisfies the conditions:
there exists such neighborhood of the point $(a=0,b=0)$ that
$f_{aa}(a,0)\not\equiv 0$ and the curve determined by the equation
$b+\varepsilon f(a,b)=0$ is monotone for every sufficiently small
$\varepsilon>0$.

Then the formulas (\ref{141}) and (\ref{18}) determine different surfaces.
\end{theorem}

{\sl PROOF.} It is easy to infer from (\ref{18}) that as $\varepsilon>0$
is sufficiently small in the neighborhood of the point $(t=0,x_1=0,x_2=0)$
the singularities surface has the form
\begin{equation}
y=\frac{t}{2}-\varepsilon f(x,-\frac{t}{2})+o(\varepsilon)\ .\label{20}
\end{equation}

Now let us find the equation of singularities surface on the basis of
formulas (\ref{141}).
As $\varepsilon=0$ the solution of (\ref{141}) looks as follows
$$x(t,l)=a_-(t,l)=a_+(t,l)=l\ ;\ y(t,l)=b_-(t,l)=-b_+(t,l)=t/2\ .$$

As $\varepsilon>0$ let us seek the solution of the system (\ref{141}) in the
form
\begin{equation}
x(t,l)=l+\varepsilon \widetilde x(t,l)+o(\varepsilon)\ ;\
y(t,l)=t/2+\varepsilon \widetilde y(t,l)+o(\varepsilon)\ .
\label{21}
\end{equation}
Then taking into account the relations
$$ \begin{aligned}
\ &x(t,l)=a_-(t,l)=a_+(t,l)+t\varepsilon f_a(a_+,b_+) \\
\ &y(t,l)=b_-(t,l)=b_+(t,l)+t\left(1+\varepsilon f_b(a_+,b_+)\right) \ ,
\end{aligned} $$
one obtains $\widetilde {P_t}$, $\widetilde {I_t}$,
$\widetilde {J_t}$ and finds
\begin{equation}\begin{aligned}
\ &\frac{\partial\widetilde x(t,l)}{\partial t}=
\frac{1}{2t}\int\limits_0^t f_ad\tau
\\
\ &\widetilde y(t,l)=-f(l,0)-\frac{1}{4}\int\limits_0^t\tau f_{aa}d\tau+
\frac{1}{2}\int\limits_0^t f_bd\tau+
\frac{1}{4t}\int_0^t\tau^2 f_{aa}d\tau\ ,
\end{aligned}\label{22}
\end{equation}
where the values of the derivatives of the function $f(a,b)$ are taken at
the point $(l,-\tau/2)$. But functions (\ref{21}) must satisfy (\ref{20})
which in view of (\ref{22}) is equivalent to the following identity
$$
f(l,-t/2)-f(l,0)-\frac{1}{4}\int\limits_0^t\left(\tau-\frac{\tau^2}{t}
\right)f_{aa}d\tau+\frac{1}{2}\int\limits_0^t f_bd\tau\equiv 0\ ,$$
where the values of the derivatives of the function $f(a,b)$ are taken at
the point $(l,-\tau/2)$. Further, multiplying by $t$ then three times
differentiating with respect to $t$ and substituting $t=0$ one obtains
$f_{aa}(l,0)\equiv 0$. The last identity contradicts to the conditions on
function $f(a,b)$.

\epr

Differentiating (\ref{18}) with respect to $t$ and changing notations from
$t$ to $\tau$ one gets the following condition on the surface $\Fc$
\begin{equation}
(a_+-a_-)\left(x^v\right)_\tau+(b_+-b_-)\left(y^v\right)_\tau=
S_0(a_+,b_+)-S_0(a_-,b_-)\quad ,\label{23}
\end{equation}

%
%  THEOREM 3.2
%
\begin{theorem}
Suppose one has the potential initial data (\ref{11}). Then the surfaces
which are defined by (\ref{141}) and (\ref{23}) coincide iff the following
relation is true
\begin{equation}
\left(a_+-a_-\right)\frac{d}{dl}\int\limits_0^ti^* d\tau +
\left(b_+-b_-\right)\frac{d}{dl}\int\limits_0^tj^* d\tau =
\left(S_+-S_-\right)\frac{d}{dl}\int\limits_0^tp^* d\tau\ ,\label{26}
\end{equation}
where
$$\begin{aligned}
\ &p^*\equiv
\left(b_+-b_-\right)\frac{(a_+)_\tau+(a_-)_\tau}{2}-
\left(a_+-a_-\right)\frac{(b_+)_\tau+(b_-)_\tau}{2}
\\
\ &i^*\equiv
\left(b_+-b_-\right)\frac{(S_+)_\tau+(S_-)_\tau}{2}-
\left(S_+-S_-\right)\frac{(b_+)_\tau+(b_-)_\tau}{2}
\\
\ &j^*\equiv
\left(a_+-a_-\right)\frac{(S_+)_\tau+(S_-)_\tau}{2}-
\left(S_+-S_-\right)\frac{(a_+)_\tau+(a_-)_\tau}{2}\ .
\end{aligned}$$
\end{theorem}

{\sl PROOF.} It is easy to see that because of our construction of the
shock surface $x^s,y^s$, keeping in mind the assumption that
$x^s=x^v$, $y^s=y^v$, one has
\begin{equation}\begin{aligned}
\ &x^s=a_-+\tau u_-(a_-,b_-)=
a_++\tau u_+(a_+,b_+) \\
\ &y^s=a_-+\tau v_-(a_-,b_-)=
a_++\tau v_+(a_+,b_+)\ .
\end{aligned}\label{24}
\end{equation}
and the relations (\ref{141}) at least locally determine the shock surface
(if it exists). So (\ref{23}) occurs to be an additional relation and it
is consistent with (\ref{141}) only if it follows from them.

Hence the condition (\ref{23}) can be rewritten as follows
\begin{equation}
(a_+-a_-){\widetilde I}_t+
(b_+-b_-){\widetilde J}_t=
\left(S_0(a_+,b_+)-S_0(a_-,b_-)\right){\widetilde P}_t\quad ,\label{25}
\end{equation}
where ${\widetilde P}_t$, ${\widetilde I}_t$,
${\widetilde J}_t$ are taken from (\ref{15}). Further,
it is easy to check that the following relations are true
$$\begin{aligned}
\ &\left((a_+)_\tau (b_+)_l-(b_+)_\tau (a_+)_l\right)-
\left((a_-)_\tau (b_-)_l-(b_-)_\tau (a_-)_l\right)=
\\
\ &\frac{d}{dl}
\left[\left(b_+-b_-\right)\frac{(a_+)_\tau+(a_-)_\tau}{2}-
\left(a_+-a_-\right)\frac{(b_+)_\tau+(b_-)_\tau}{2}\right]-
\\
\ &\frac{d}{d\tau}
\left[\left(b_+-b_-\right)\frac{(a_+)_l+(a_-)_l}{2}-
\left(a_+-a_-\right)\frac{(b_+)_l+(b_-)_l}{2}\right]=
\frac{d}{dl}p^*-\frac{d}{d\tau}p\ ;
\end{aligned}$$

$$\begin{aligned}
\ &u_+\left((a_+)_\tau (b_+)_l-(b_+)_\tau (a_+)_l\right)-
u_-\left((a_-)_\tau (b_-)_l-(b_-)_\tau (a_-)_l\right)=
\\
\ &\left((S_+)_\tau (b_+)_l-(b_+)_\tau (S_+)_l\right)-
\left((S_-)_\tau (b_-)_l-(b_-)_\tau (S_-)_l\right)=
\\
\ &\frac{d}{dl}
\left[\left(b_+-b_-\right)\frac{(S_+)_\tau+(S_-)_\tau}{2}-
\left(S_+-S_-\right)\frac{(b_+)_\tau+(b_-)_\tau}{2}\right]-
\\
\ &\frac{d}{d\tau}
\left[\left(b_+-b_-\right)\frac{(S_+)_l+(S_-)_l}{2}-
\left(S_+-S_-\right)\frac{(b_+)_l+(b_-)_l}{2}\right]=
\frac{d}{dl}i^*-\frac{d}{d\tau}i\ ;
\end{aligned}$$

$$\begin{aligned}
\ &-\left[v_+\left((a_+)_\tau (b_+)_l-(b_+)_\tau (a_+)_l\right)-
v_-\left((a_-)_\tau (b_-)_l-(b_-)_\tau (a_-)_l\right)\right]=
\\
\ &\left((S_+)_\tau (a_+)_l-(a_+)_\tau (S_+)_l\right)-
\left((S_-)_\tau (a_-)_l-(a_-)_\tau (S_-)_l\right)=
\\
\ &\frac{d}{dl}
\left[\left(a_+-a_-\right)\frac{(S_+)_\tau+(S_-)_\tau}{2}-
\left(S_+-S_-\right)\frac{(a_+)_\tau+(a_-)_\tau}{2}\right]-
\\
\ &\frac{d}{d\tau}
\left[\left(a_+-a_-\right)\frac{(S_+)_l+(S_-)_l}{2}-
\left(S_+-S_-\right)\frac{(a_+)_l+(a_-)_l}{2}\right]=
\frac{d}{dl}j^*-\frac{d}{d\tau}j\ ,
\end{aligned}$$
where $S_k\equiv S_0(a_k,b_k)$, $(k=+,-)$. Taking into account these
formulas relation (\ref{25}) produce (\ref{26}) and inversely.

\epr

%
%  REMARK 3.1
%
\begin{remark}
In case the time interval $t$ is small and initial data are sufficiently
smooth all terms in relation (\ref{26}) approximately equal to zero and
one recovers formulas suggested in \cite{RykoUMNa96}. Thus when the time
elapsing from the moment of the emerging of a shock front is small the
variational representation is approximately valid.
\end{remark}

%
%  REMARK 3.2
%
\begin{remark}
The nontrivial 2-D example which is known to the author where the
variational representation works rigorously is the 2-D Riemann
problem for (\ref{1}) with compressive piecewise constant potential
velocity vector. The example of corresponding initial potential $S_0(a,b)$
is shown below ($A,B,C$ are constants)
\begin{align*}
\ & S_0(a,b)=Ca-\left(B+\frac{A+C}{\sqrt{3}}\right)b\ ;\
b\ge\max\left(0,\sqrt{3}a\right) \\
\ & S_0(a,b)=-Aa-Bb\ ;\
a\ge 0,\sqrt{3}a\ge b\ge -\sqrt{3}a \\
\ & S_0(a,b)=Ca+\left(\frac{A+C}{\sqrt{3}}-B\right)b\ ;\
b\le\min\left(0,-\sqrt{3}a\right)\ ,\\
\ & C>0, A+\sqrt{3}B>0, A-\sqrt{3}B>0\ .
\end{align*}
This result is obtained when the discontinuity lines are
straight and between every two neighboring lines is the same angle:
$90^\circ$ or $120^\circ$ or $180^\circ$. The conjecture is that the angle
value can be taken arbitrary.

But there the another type of singularity arises: $\delta$-function in one
point for density.

The validity of variational principle in this case could be understood in
the following way: the geometry of corresponding initial potential was
flat and had high degree of symmetry.
\end{remark}

%
%  SECTION 4
%
\section{The flow description inside the shock.}

%
%  THEOREM 4.1
%
\begin{theorem}
Suppose there exists the generalized solution to the problem (\ref{1}),
(\ref{11}) in the form (\ref{12}), (\ref{141}), (\ref{15}) and the surface
$\Sc$ can be defined from the parametric equations $x=x(t,l)$,
$y=y(t,l)$.  Then the following system of equations is true
\begin{align}
\ &{\left(\widetilde {P_t}\right)}^\cdot+
x_l\left\{V (\varrho_+-\varrho_-)-(\varrho_+v_+-\varrho_-v_-)\right\}+
\notag \\
\ &y_l\left\{(\varrho_+u_+-\varrho_-u_-)-U(\varrho_+-\varrho_-)\right\}=0
\notag \\ \notag \\
\ &{\left(\widetilde {I_t}\right)}^\cdot+ x_l\left\{
V(\varrho_+u_+-\varrho_-u_-)-(\varrho_+u_+v_+-\varrho_-u_-v_-)\right\}+
\notag \\
\ &y_l\left\{(\varrho_+u_+^2-\varrho_-u_-^2)-
U(\varrho_+u_+-\varrho_-u_-)\right\}=0
\label{27} \\ \notag \\
\ &{\left(\widetilde {J_t}\right)}^\cdot+ x_l\left\{
V(\varrho_+v_+-\varrho_-v_-)-(\varrho_+v_+^2-\varrho_-v_-^2)\right\}+
\notag \\
\ &y_l\left\{(\varrho_+u_+v_+-\varrho_-u_-v_-)-
U(\varrho_+v_+-\varrho_-v_-)\right\}=0
\notag \\ \notag \\
\ &\dot y=V\quad ,\quad \dot x=U\ , \notag
\end{align}
where $U=\widetilde {I_t}/\widetilde {P_t}$,
$V=\widetilde {J_t}/\widetilde {P_t}$ and 'dot' denotes the
differentiation with respect to $t$.
\end{theorem}

{\sl PROOF.} Taking into account (\ref{4}), (\ref{5}) one has
\begin{equation}\begin{aligned}
\ &x(t,l)=a_i+tu_0^i(a_i,b_i) \\
\ &y(t,l)=b_i+tv_0^i(a_i,b_i) \\
\ &\varrho_i=\varrho_0(a_i,b_i)/D(a_i,b_i)\ ,
\end{aligned}\label{28}
\end{equation}
where $(i=+,-)$ and
$$
D(a_i,b_i)\equiv \left(1+t\left(u_0^i\right)_a\right)
\left(1+t\left(v_0^i\right)_b\right)-t^2
\left(u_0^i\right)_b\left(v_0^i\right)_a\ .
$$

Further in this proof for our convenience we will do all calculations
independently for indices "$+$", "$-$" and so omit indices in the
expressions and write $a$, $b$, $\varrho_0$, $u$, $v$, $D$ instead of
$a_i$, $b_i$, $\varrho_0^i$, $u_0^i$, $v_0^i$, $D_i$, $i=+,-$.
From (\ref{28}) one can infer after differentiation the
relations for $x_l$, $y_l$, $\dot x$, $\dot y$.

Now consider the right hand side of the first equation from (\ref{27})
and write all terms concerning the index "$+$" (for "$-$" the calculations
are analogous). One has
$$\begin{aligned}
\ &\varrho\left[x_l\left(V-v\right)+y_l\left(u-U\right)\right]=
\varrho\left[x_l\left(\dot y-v\right)+y_l\left(u-\dot x\right)\right]=
\\
\ &\varrho\left[a_l\left(1+tu_a\right)+b_ltu_b\right]
\left[\dot b\left(1+tv_b\right)+\dot atv_a\right]-
\\
\ &\varrho\left[a_ltv_a+b_l\left(1+tv_b\right)\right]
\left[\dot a\left(1+tu_a\right)+\dot btu_b\right]=
\\
\ &-\varrho D\left(\dot ab_l-\dot ba_l\right)=
-\varrho_0\left(\dot ab_l-\dot ba_l\right)\ .
\end{aligned}$$
Thus the first equation from (\ref{27}) is equivalent to the first
equation from (\ref{15}). Taking into account the relations
$$\begin{aligned}
\ &\varrho u\left[x_l\left(V-v\right)+y_l\left(u-U\right)\right]=
-\varrho_0u \left(\dot ab_l-\dot ba_l\right)
\end{aligned}$$
and
$$\begin{aligned}
\ &\varrho v\left[x_l\left(V-v\right)+y_l\left(u-U\right)\right]=
-\varrho_0v \left(\dot ab_l-\dot ba_l\right)
\end{aligned}$$
one obtains the other two equations.

\epr

\begin{corollary}
The system (\ref{27}) is well defined.
\end{corollary}

{\sl PROOF.} We have to prove that $\widetilde {P_t} >0$ under the
dynamics of (\ref{27}). Let us integrate $\widetilde {P_t}$ with respect to
$l$ from some $l_0$ to $l_0+\Delta l$, where $\Delta l$ is small enough.
Note that $\int_{l_0}^{l_0+\Delta l}\widetilde {P_t}dl$ is exactly the sum
of the areas which are bounded by the curves ($i=+,-$):
$$\left(a_i(\tau,l_0),b_i(\tau,l_0)\right)\ ,\
\left(a_i(\tau,l_0+\Delta l),b_i(\tau,l_0+\Delta l)\right)\ ,\quad
0\le\tau\le t\ ;$$
$$\left(a_i(t,l),b_i(t,l)\right)\ ,\
\left(a_i(0,l),b_i(0,l)\right)\ ,\quad
l_0\le l\le l_0+\Delta l\ .$$
Since $\Delta l$ is arbitrary then $\widetilde {P_t}>0$ as $t>0$.

\epr

%
%  REMARK 4.1
%
\begin{remark}
The system (\ref{27}) is nonstrictly hyperbolic system with one
eigenvalue and three eigenvectors. From the first three equations of
system (\ref{27}) it is obvious to find the additional relation between
${\left(\widetilde {P_t}\right)}^\cdot$, ${\left(\widetilde
{I_t}\right)}^\cdot$, ${\left(\widetilde {J_t}\right)}^\cdot$ by
eliminating $x_l$ and $y_l$.
\end{remark}

%
%  REMARK 4.2
%
\begin{remark}
From the Cauchy-Kovalevskaya theorem one immediately gets the local
existence and uniqueness theorem for (\ref{27}) in case of analytic
coefficients and initial data.
\end{remark}

%
%  SECTION 5
%
\section{Specific type of motion. Constant external state.}

The system (\ref{27}) does not satisfy Friedrichs' symmetrizability
condition so to obtain the existence theorem encounters some problems.
Nevertheless the internal dynamics of (\ref{27}) is highly nontrivial
which we demonstrate here in the simplest case of constant external
density and velocity.

Suppose that in (\ref{27}) $\varrho_+\equiv const\equiv\widetilde\varrho$,
$\varrho_-\equiv const\equiv\varrho$,
$u_+=v_+\equiv 0$, $u_-\equiv const\equiv u$, $v_-\equiv const\equiv v$;
the velocity vector $(u,v)$ satisfies condition II).
To simplify the notations we also drop index '$t$' and 'waves' in
(\ref{27}). Then one has
\begin{equation}\begin{aligned}
\ &\dot P+\left(\varrho-\widetilde\varrho\right)\left\{y_lU-x_lV\right\}-
\varrho\left\{y_lu-x_lv\right\}=0 \\
\ &\dot I+\varrho u\left\{y_lU-x_lV\right\}-
\varrho u\left\{y_lu-x_lv\right\}=0 \\
\ &\dot J+\varrho v\left\{y_lU-x_lV\right\}-
\varrho v\left\{y_lu-x_lv\right\}=0 \\
\ &\dot x=U\equiv I/P\quad ,\quad\dot y=V\equiv J/P\ .
\end{aligned}\label{29}
\end{equation}
From (\ref{29}) one immediately gets the first integrals
\begin{equation}
uJ-vI=uJ_0-vI_0\equiv C(l)\ .\label{30}
\end{equation}
and
\begin{equation}
u\dot y-v\dot x=\frac{C(l)}{P}\ .\label{31}
\end{equation}

Now let us take the special initial conditions such that $C(l)\equiv 0$.
In other words there exists such $k_0(l)$ that
\begin{equation}
I_0=k_0u\quad ,\quad J_0=k_0v\ .\label{32}
\end{equation}
Then taking into account (\ref{30}) there exists some unknown $k(t,l)$
such that
$$
I=ku\quad ,\quad J=kv\ .
$$
And from (\ref{31})
$$
uy_l-vx_l=uy_0^\prime-vx_0^\prime\equiv G(l)>0\ .
$$
due to condition II).
After rather simple transformations one arrives to the following system
\begin{equation}\begin{aligned}
\ &\dot P+G(l)\left\{\widehat k\left(\varrho-\widetilde\varrho\right)-
\varrho\right\}=0 \\
\ &\left(P\widehat k\right)^\cdot+\varrho G(l)\left\{\widehat
k-1\right\}=0\ ,
\end{aligned}\label{33}
\end{equation}
where $\widehat k\equiv k/P$.

Now expressing $\widehat k$ from the first equation of (\ref{33}) and
substituting it into the second equation one can find the expressions for
$P$ and $\widehat k$ which read
\begin{equation}\begin{aligned}
\ &P^2=P_0(l)^2-2G(l)P_0(l)N(l)t+\varrho\widetilde\varrho G(l)^2t^2 \\
\ &\widehat k\left(\varrho-\widetilde\varrho\right)=
\varrho-\frac{\dot P}{G(l)}\ ,
\end{aligned}\label{34}
\end{equation}
where $N(l)\equiv {\widehat
k}_0(l)\left(\varrho-\widetilde\varrho\right)-\varrho$.

From (\ref{29}) one obtains that $\dot x=\widehat ku$, $\dot y=\widehat
kv$ and the stability condition (see Definition 2) takes the form
\begin{equation}
0<\int\limits_0^t\widehat k(\tau,l)d\tau<t\ .
\label{35}\end{equation}

It is easy to find that from (\ref{35}) taking into account (\ref{34})
follows the stability condition (see Definition 2)
\begin{equation}
0<{\widehat k}_0(l)<1\ .
\label{36}\end{equation}
Then from (\ref{36}) one can easily see that $N(l)<0$ and so (\ref{34})
are well defined.  From (\ref{34}) also follows that $\widehat k\to
\varkappa$ as $t\to\infty$, where
\begin{equation}
\varkappa\equiv\frac{\sqrt{\varrho}}{\sqrt{\varrho}+
\sqrt{\widetilde\varrho}}\ .
\label{37}\end{equation}

Thus we have proved the following theorem

%
%  THEOREM 5.1
%
\begin{theorem}
Suppose that $\varrho_+\equiv const\equiv\widetilde\varrho$,
$\varrho_-\equiv const\equiv\varrho$,
$u_+=v_+\equiv 0$, $u_-\equiv const\equiv u$, $v_-\equiv const\equiv v$.

Suppose also that $G(l)>0$ and (\ref{36}) is true. Than there exists the
solution to the problem (\ref{1}), (\ref{11}) and the shock front tends to
the following one as $t\to +\infty$
$$\begin{aligned}
\ &x(l,t)=x_0(l)+\varkappa ut \\
\ &y(l,t)=y_0(l)+\varkappa vt\ ,
\end{aligned}$$
where $\varkappa$ is taken from (\ref{37}).
\end{theorem}

Finally to illustrate the nontrivial character of the problem even with
constant external fields let us derive the equation in the case
$C(l)\not\equiv 0$, but $\varrho=\widetilde\varrho$. Then from (\ref{29})
one infers
\begin{equation}\begin{aligned}
\ &\dot P=\varrho\left\{y_lu-x_lv\right\} \\
\ &\dot I+\dot P\left(\frac{I}{P}-u\right)=\varrho x_l\frac{C(l)}{P} \\
\ &\dot x=\frac{I}{P}\ .
\end{aligned}\label{38}
\end{equation}
Differentiating the first equation from (\ref{38}) with respect to $t$ and
taking into account integral (\ref{31}) one gets the equation for $P$
\begin{equation}
\ddot P=\varrho\left(\frac{C(l)}{P}\right)_l
\label{39}\end{equation}
So even in the simplest case our system delivers us rather unusual
equations of type (\ref{39}). To illustrate this let us perform the
stability analysis for small perturbations to the model linear equation
\begin{equation}
\ddot P=KP_x\ ,
\label{40}\end{equation}
where $K=const$. One have to find partial solutions to (\ref{40}) in the
form
\begin{equation}
P(t,x)=e^{i(\xi x+\lambda t)}\ ,
\label{41}\end{equation}
where $i^2=-1$, $\xi\in \R$, $\lambda=\sigma+i\Delta$, $\sigma\in \R$,
$\Delta\in \R$. One immediately gets the restriction to the choice of
$\xi$ and $\lambda$
$$
\lambda^2=-Ki\xi\ .
$$
Further, one has
\begin{equation}
\sigma^2=\Delta^2\quad ;\quad 2\sigma\Delta=-K\xi\ .
\label{42}\end{equation}
From (\ref{42}) it follows that
$$
2\Delta^2=\pm K\xi\ ,
$$
so we can choose such signs that $-\Delta\sim const\sqrt{\xi}$ and
$-\Delta$ tends to $+\infty$ as $\xi$ tends to $+\infty$. Thus in
(\ref{41}) one has an arbitrary rapid growth of small perturbations with
high frequencies and equation (\ref{40}) is ill-posed in the class of
functions of finite smoothness because it does not satisfy classical
Petrovsky condition.

%
%  APPENDIX
%
\section*{Appendix.}
\addcontentsline{toc}{section}{\bf Appendix.}

Let us carry out some heuristic calculations to obtain
the system (\ref{27}). The generalized solution (\ref{1}), (\ref{11}) can
also be written in the form
\begin{equation}\begin{aligned}
\ &u=u_-(t,x,y)+\left(u_+(t,x,y)-u_-(t,x,y)
\right)H(S) \\
\ &v=v_-(t,x,y)+\left(v_+(t,x,y)-v_-(t,x,y)
\right)K(S) \\
\ &\varrho=\varrho_-(t,x,y)+\left(\varrho_+(t,x,y)-\varrho_-(t,x,y)
\right)R(S)+\lambda|_{S=0}\delta(S)\ ,
\end{aligned}\label{43}
\end{equation}
where $\delta$ is usual Dirac $\delta$-function, but $H$, $K$, $R$ are
different Heaviside functions which can be distinguished by means of the
following heuristic multiplication formulas
\begin{equation}
\delta\cdot H=s_1\cdot\delta\quad ,\quad \delta\cdot K=s_2\cdot\delta\ ,
\label{44}
\end{equation}
where $s_1,s_2$ are some functions on the surface ${\cal S}$ (in what
follows we need not multiplication with $R$, so it is not included in
(\ref{44})). In addition the following rather natural formulas are
supposed to be true in the sense of distributions
\begin{equation}
H^2\approx K^2\approx HK\approx RH\approx RK\approx H\ .\label{45}
\end{equation}
Let us note that these formulas can be treated rigorously, for example,
with the help of theory of new generalized functions \cite{ColoBElI85},
\cite{BiagBNTh90} (all basic ideas and lines in application to physics
can also be found in \cite{ColoBAMS90}) but here we do not need such
rigor.

We can write
$$
U\equiv u_-+s_1(u_+-u_-)\quad ,\quad V\equiv
v_-+s_2(v_+-v_-)\ .
$$

Then  taking into account the relations (\ref{45}) one has the following
equalities
\begin{align*}
\ &\varrho u=\varrho_-u_-+(\varrho_+u_+-\varrho_-u_-)H+
\lambda U\delta \\
\ &\varrho v=\varrho_-v_-+(\varrho_+v_+-\varrho_-v_-)H+
\lambda V\delta
\\ \\
\ &\varrho u^2=\varrho_-u_-^2+(\varrho_+u_+^2-\varrho_-u_-^2)H+
\lambda U^2\delta \\
\ &\varrho uv=\varrho_-u_-v_-+
(\varrho_+u_+v_+-\varrho_-u_-v_-)H+
\lambda UV\delta \\
\ &\varrho v^2=\varrho_-v_-^2+(\varrho_+v_+^2-\varrho_-v_-^2)H+
\lambda V^2\delta
\end{align*}
Substituting these equalities in the system (\ref{1}) and equating to zero
expressions with different kind of singularities one obtains
\begin{equation}\begin{aligned}
\ &S_t+US_x+VS_y=0 \\ \\
\ &\lambda_t+\left(\lambda U\right)_x+\left(\lambda V\right)_y+
(\varrho_+-\varrho_-)S_t+ \\
\ &(\varrho_+u_+-\varrho_-u_-)S_x+
(\varrho_+v_+-\varrho_-v_-)S_y=0 \\ \\
\ &\left(\lambda U\right)_t+\left(\lambda U^2\right)_x+\left(\lambda
UV\right)_y+(\varrho_+u_+-\varrho_-u_-)S_t+ \\
\ &(\varrho_+u_+^2-\varrho_-u_-^2)S_x+
(\varrho_+u_+v_+-\varrho_-u_-v_-)S_y=0 \\ \\
\ &\left(\lambda V\right)_t+\left(\lambda UV\right)_x+
\left(\lambda V^2\right)_y+(\varrho_+v_+-\varrho_-v_-)S_t+
\\
\ &(\varrho_+u_+v_+-\varrho_-u_-v_-)S_x+
(\varrho_+v_+^2-\varrho_-v_-^2)S_y=0\ .
\end{aligned}\label{46}
\end{equation}

Now let us mention that the surface $\Sc$ can be defined from the
equation $S\equiv x-X(t,y)=0$ but the functions $\lambda$, $U$, $V$ depend
only on $t,y$. Introducing the differentiation along the direction
$\left(U,V\right)$ which will be denoted by 'dot' rewrite the system
(\ref{46}) in the form
\begin{align}
\ &\dot y=V\quad ,\quad \dot X=U\quad ,\quad y(0,l)=l \notag \\
\ &\dot\lambda+\lambda V_y+X_y\left(V (\varrho_+-\varrho_-)-
(\varrho_+v_+-\varrho_-v_-)\right)+ \notag \\
\ &(\varrho_+u_+-\varrho_-u_-)-U(\varrho_+-\varrho_-)=0 \notag \\
\notag \\
\ &\lambda\dot {U}+X_y\bigl[V\left(\varrho_+u_+-\varrho_-u_-
-U(\varrho_+-\varrho_-)\right)+ \notag \\
\ &U(\varrho_+v_+-\varrho_-v_-)-
(\varrho_+u_+v_+-\varrho_-u_-v_-)\bigr]+ \notag \\
\ &U\left(U(\varrho_+-\varrho_-)-(\varrho_+u_+-\varrho_-u_-)
\right)+ \notag \\
\ &(\varrho_+u_+^2-\varrho_-u_-^2)-
U(\varrho_+u_+-\varrho_-u_-)=0 \label{47}\\ \notag \\
\ &\lambda\dot {V}+X_y\bigl[V\left(\varrho_+v_+-\varrho_-v_-
-V(\varrho_+-\varrho_-)\right)+ \notag \\
\ &V(\varrho_+v_+-\varrho_-v_-)-
(\varrho_+v_+^2-\varrho_-v_-^2)\bigr]+ \notag \\
\ &V\left(U(\varrho_+-\varrho_-)-(\varrho_+u_+-\varrho_-u_-)
\right)+ \notag \\
\ &(\varrho_+u_+v_+-\varrho_-u_-v_-)-
U(\varrho_+v_+-\varrho_-v_-)=0\ . \notag
\end{align}

Now let us transform the system (\ref{47}) in the following way. Multiply
the first equation by $U$ and add to the second equation, then multiply
the first equation by $V$ and add to the third equation. Finally
multiplying all obtained equations by $y_l$ and using the relation $\dot
{y_l}=V_yy_l$ one gets the system (\ref{27}).

{\bf Acknowledgments.} The author is grateful to E.~Aurell, C.~Boldrighini,
E.~Caglioti, V.~Oseledets, M.~Pulverenti and Ya.~Sinai for useful
discussions and encouragement, to {\it Dipartimento di Matematica,
Universit\`a di Roma "La Sapien\-za" (Italy)} for warm hospitality. This
work was also supported by {\it Russian Foundation for Basic
Researches}, grants No. 99--01--00314, No. 00--01--00387.

%
%  BIBLIOGRAPHY
%

%
\vfill\eject
\end{document}